\documentclass{amsart}
\usepackage{amstext}
\usepackage{amsthm}
\usepackage{amssymb}

\makeatletter

\numberwithin{equation}{section}
\numberwithin{figure}{section}
\theoremstyle{plain}
\newtheorem{thm}{Theorem}

\allowdisplaybreaks
\usepackage[sort&compress,square,numbers]{natbib}

\makeatother

\global\long\def\relphantom#1{\mathrel{\phantom{{#1}}}}

\usepackage{tikz}
\usetikzlibrary{matrix,positioning}

\begin{document}

\title[Nonlinear differential equations and Boole numbers]{Nonlinear differential equations arising from Boole numbers and their applications}

\author{Taekyun Kim}
\address{Department of Mathematics, Kwangwoon University, Seoul 139-701, Republic
	of Korea}
\email{tkkim@kw.ac.kr}

\author{Dae San Kim}
\address{Department of Mathematics, Sogang University, Seoul 121-742, Republic
	of Korea}
\email{dskim@sogang.ac.kr}
\begin{abstract}
In this paper, we study nonlinear differential equations satisfied by the generating function of Boole numbers.
In addition, we derive some explicit and new interesting identities involving 
Boole numbers and higher-order Boole numbers arising from our nonlinear differential equations.
\end{abstract}



\maketitle

\section{Introduction}

The Boole polynomials, $Bl_{n}\left(x\mid\lambda\right)$,
$\left(n\ge0\right)$, are given by the generating function 
\begin{equation}
\frac{1}{1+\left(1+t\right)^{\lambda}}\left(1+t\right)^{x}=\sum_{n=0}^{\infty}Bl_{n}\left(x\mid\lambda\right)\frac{t^{n}}{n!},\quad\left(\text{see \cite{key-5,key-6,key-7,key-8,key-9,key-10,key-18}}\right),\label{eq:1}
\end{equation}
where we assume that $\lambda\neq 0$.

When $x=0$, $Bl_{n}\left(x\right)=Bl_{n}\left(0\mid\lambda\right)$,
$\left(n\ge0\right)$, are called the Boole numbers. The higher-order
Boole polynomials (or Peters polynomials) are also defined by the
generating function 
\begin{equation}
\left(\frac{1}{1+\left(1+t\right)^{\lambda}}\right)^{r}\left(1+t\right)^{x}=\sum_{n=0}^{\infty}Bl_{n}^{\left(r\right)}\left(x\mid\lambda\right)\frac{t^{n}}{n!},\quad\left(r\in\mathbb{N}\right),\quad\left(\text{see \cite{key-18} }\right).\label{eq:2}
\end{equation}

The first few Boole and higher-order Boole polynomials are as follows:
\[
Bl_{0}\left(x\mid\lambda\right)=\frac{1}{2},\quad Bl_{1}\left(x\mid\lambda\right)=\frac{1}{4}\left(2x-\lambda\right),\quad Bl_{2}\left(x\mid\lambda\right)=\frac{1}{4}\left(2x\left(x-\lambda-1\right)+\lambda\right),
\]
and 
\begin{align*}
Bl_{0}^{\left(r\right)}\left(x\mid\lambda\right) & =2^{-r},\quad Bl_{1}^{\left(r\right)}\left(x\mid\lambda\right)=2^{-\left(r+1\right)}\left(2x-\lambda\right),\\
Bl_{2}^{\left(r\right)}\left(x\mid\lambda\right) & =2^{-\left(r+2\right)}\left(4x\left(x-1\right)+\left(2-4x\right)\lambda r+r\left(r-1\right)\lambda^{2}\right),\cdots.
\end{align*}

With the viewpoint of umbral calculus, Boole numbers and polynomials
have been studied by several authors (see \cite{key-1,key-2,key-3,key-4,key-5,key-6,key-7,key-8,key-9,key-10,key-11,key-12,key-13,key-14,key-15,key-16,key-17,key-18,key-19,key-20}). 

Recently, Kim-Kim has studied the following nonlinear differential
equations(see \cite{key-7,key-9}): 

\begin{equation}
\left(\frac{d}{dt}\right)^{N}F\left(t\right)=\frac{\left(-1\right)^{N}}{\left(1+t\right)^{N}}\sum_{j=2}^{N+1}\left(j-1\right)!\left(N-1\right)!H_{N-1,j-2}F\left(t\right)^{j},\quad\left(N\in\mathbb{N}\right),\label{eq:3}
\end{equation}
where 
\begin{align*}
H_{N,0} & =1,\quad\text{for all }N,\\
H_{N,1} & =H_{N}=1+\frac{1}{2}+\cdots+\frac{1}{N},\\
H_{N,j} & =\frac{H_{N-1,j-1}}{N}+\frac{H_{N-2,j-1}}{N-1}+\cdots+\frac{H_{0,j-1}}{1},\quad H_{0,j-1}=0,\quad\left(2\le j\le N\right).
\end{align*}

From (\ref{eq:3}), they derived some explicit and new identities for the
Bernoulli numbers of the second kind and the higher-order Bernoulli numbers of the second kind. 

The purpose of this paper is to give some explicit and new identities for
the Boole numbers and the higher-order Boole numbers arising from nonlinear differential equations.

\section{Nonlinear differential equations arising from the generating function of Boole numbers}

Let 
\begin{equation}
F=F\left(t ;\lambda\right)=\frac{1}{\left(1+t\right)^{\lambda}+1}.\label{eq:4}
\end{equation}

Then, by (\ref{eq:4}), we get 
\begin{align}
F^{\left(1\right)} & =\frac{d}{dt}F\left(t\right)\label{eq:5}\\
 & =\left(\frac{1}{\left(1+t\right)^{\lambda}+1}\right)^{2}\frac{\left(-1\right)\lambda}{\left(1+t\right)}\left(1+t\right)^{\lambda}\nonumber \\
 & =\frac{\left(-1\right)\lambda}{1+t}\frac{1}{\left(\left(1+t\right)^{\lambda}+1\right)^{2}}\left(\left(1+t\right)^{\lambda}-1+1\right)\nonumber \\
 & =\frac{\left(-1\right)\lambda}{1+t}\left(F-F^{2}\right),\nonumber 
\end{align}
and 
\begin{align}
F^{\left(2\right)} & =\frac{dF^{\left(1\right)}}{dt}\label{eq:6}\\
 & =\frac{\left(-1\right)^{2}\lambda}{\left(1+t\right)^{2}}\left(F-F^{2}\right)-\frac{\lambda}{1+t}\left(F^{\left(1\right)}-2FF^{\left(1\right)}\right)\nonumber \\
 & =\frac{\left(-1\right)^{2}\lambda}{\left(1+t\right)^{2}}\left(F-F^{2}\right)+\frac{\left(-1\right)^{2}\lambda}{\left(1+t\right)^{2}}\left(1-2F\right)\left(F-F^{2}\right)\nonumber \\
 & =\frac{\left(-1\right)^{2}\lambda}{\left(1+t\right)^{2}}\left\{ \left(1+\lambda\right)F-\left(1+3\lambda\right)F^{2}+2\lambda F^{3}\right\} .\nonumber 
\end{align}

Continuing this process, we set 
\begin{equation}
F^{\left(N\right)}=\left(\frac{d}{dt}\right)^{N}F\left(t\right)=\frac{\left(-1\right)^{N}\lambda}{\left(1+t\right)^{N}}\sum_{i=1}^{N+1}a_{i-1}\left(N;\lambda\right)F^{i},\label{eq:7}
\end{equation}
where $N=0,1,2,\dots$. 

From (\ref{eq:7}), we have 
\begin{align}
 & \relphantom =F^{\left(N+1\right)}\label{eq:8}\\
 & =\frac{d}{dt}F^{\left(N\right)}\nonumber \\
 & =\frac{\left(-1\right)^{N+1}\lambda N}{\left(1+t\right)^{N+1}}\sum_{i=1}^{N+1}a_{i-1}\left(N;\lambda\right)F^{i}+\frac{\left(-1\right)^{N}\lambda}{\left(1+t\right)^{N}}\sum_{i=1}^{N+1}a_{i-1}\left(N;\lambda\right)iF^{i-1}F^{\left(1\right)}\nonumber \\
 & =\frac{\left(-1\right)^{N+1}\lambda N}{\left(1+t\right)^{N+1}}\sum_{i=1}^{N+1}a_{i-1}\left(N;\lambda\right)F^{i}+\frac{\left(-1\right)^{N+1}\lambda^{2}}{\left(1+t\right)^{N+1}}\sum_{i=1}^{N+1}ia_{i-1}\left(N;\lambda\right)F^{i-1}\left(F-F^{2}\right)\nonumber \\
 & =\frac{\left(-1\right)^{N+1}\lambda}{\left(1+t\right)^{N+1}}\left\{ \sum_{i=1}^{N+1}\left(N+i\lambda\right)a_{i-1}\left(N;\lambda\right)F^{i}-\sum_{i=2}^{N+2}\left(i-1\right)\lambda a_{i-2}\left(N;\lambda\right)F^{i}\right\} \nonumber \\
 & =\frac{\left(-1\right)^{N+1}\lambda}{\left(1+t\right)^{N+1}}\left\{ \left(N+\lambda\right)a_{0}\left(N;\lambda\right)F-\left(N+1\right)\lambda a_{N}\left(N;\lambda\right)F^{N+2}\right.\nonumber \\
 & \relphantom =\left.+\sum_{i=2}^{N+1}\left(\left(N+i\lambda\right)a_{i-1}\left(N;\lambda\right)-\left(i-1\right)\lambda a_{i-2}\left(N;\lambda\right)F^{i}\right)\right\} .\nonumber 
\end{align}

On the other hand, replacing $N$ by $N+1$ in (\ref{eq:7}), we get 
\begin{equation}
F^{\left(N+1\right)}=\frac{\left(-1\right)^{N+1}\lambda}{\left(1+t\right)^{N+1}}\sum_{i=1}^{N+2}a_{i-1}\left(N+1;\lambda\right)F^{i}.\label{eq:9}
\end{equation}

From (\ref{eq:8}) and (\ref{eq:9}), we can derive the following
relations: 
\begin{align}
a_{0}\left(N+1;\lambda\right) & =\left(N+\lambda\right)a_{0}\left(N;\lambda\right),\label{eq:10}\\
a_{N+1}\left(N+1;\lambda\right) & =-\left(N+1\right)\lambda a_{N}\left(N;\lambda\right)\label{eq:11}
\end{align}
and 
\begin{equation}
a_{i-1}\left(N+1;\lambda\right)=-\left(i-1\right)\lambda a_{i-2}\left(N;\lambda\right)+\left(N+i\lambda\right)a_{i-1}\left(N;\lambda\right),\label{eq:12}
\end{equation}
where $2\le i\le N+1$. 

By (\ref{eq:4}) and (\ref{eq:7}), it is easy to show that 
\begin{equation}
F=F^{\left(0\right)}=\lambda a_{0}\left(0;\lambda\right)F.\label{eq:13}
\end{equation}

By comparing the coefficients on both sides of (\ref{eq:13}),
we have 
\begin{equation}
a_{0}\left(0;\lambda\right)=\frac{1}{\lambda}.\label{eq:14}
\end{equation}

From (\ref{eq:5}) and (\ref{eq:7}), we note that 
\begin{align}
\frac{\left(-1\right)\lambda}{1+t}\left(F-F^{2}\right) & =F^{\left(1\right)}\label{eq:15}\\
& =\frac{\left(-1\right)\lambda}{1+t}\left(a_{0}\left(1;\lambda\right)F+a_{1}\left(1;\lambda\right)F^{2}\right).\nonumber 
\end{align}

Thus, by (\ref{eq:15}), we get
\[
a_{0}\left(1;\lambda\right)=1,\text{ and }a_{1}\left(1;\lambda\right)=-1.
\]
 
\begin{align}
a_{0}\left(N+1;\lambda\right) & =\left(N+\lambda\right)a_{0}\left(N;\lambda\right)\label{eq:17}\\
 & =\left(N+\lambda\right)\left(N+\lambda-1\right)a_{0}\left(N-1;\lambda\right)\nonumber \\
 & \vdots\nonumber \\
 & =\left(N+\lambda\right)\left(N+\lambda-1\right)\cdots\left(1+\lambda\right)a_{0}\left(1;\lambda\right)\nonumber \\
 & =\left(N+\lambda\right)\left(N+\lambda-1\right)\cdots\left(1+\lambda\right)\cdot1\nonumber \\
 & =\left(N+\lambda\right)_{N},\nonumber 
\end{align}
and 
\begin{align}
a_{N+1}\left(N+1;\lambda\right) & =-\left(N+1\right)\lambda a_{N}\left(N;\lambda\right)\label{eq:18}\\
 & =\left(-1\right)^{2}\lambda^{2}\left(N+1\right)Na_{N-1}\left(N-1;\lambda\right)\nonumber \\
 & \vdots\nonumber \\
 & =\left(-1\right)^{N}\lambda^{N}\left(N+1\right)N\cdots2a_{1}\left(1;\lambda\right)\nonumber \\
 & =\left(-1\right)^{N+1}\lambda^{N}\left(N+1\right)!,\nonumber 
\end{align}
where 
\[
\left(x\right)_{n}=x\left(x-1\right)\left(x-2\right)\cdots\left(x-n+1\right),\quad\left(n\ge0\right).
\]

From (\ref{eq:12}), we can derive the following equations: 
\begin{align}
&\relphantom{=}a_{1}\left(N+1;\lambda\right) \label{eq:19}\\
& =-\lambda a_{0}\left(N;\lambda\right)+\left(N+2\lambda\right)a_{1}\left(N;\lambda\right)\nonumber\\
 & =-\lambda a_{0}\left(N;\lambda\right)+\left(N+2\lambda\right)\left\{ -\lambda a_{0}\left(N-1;\lambda\right)+\left(\left(N-1\right)+2\lambda\right)a_{1}\left(N-1;\lambda\right)\right\} \nonumber \\
 & =-\lambda\left(a_{0}\left(N;\lambda\right)+\left(N+2\lambda\right)a_{0}\left(N-1;\lambda\right)\right)+\left(N+2\lambda\right)\left(N+2\lambda-1\right)a_{1}\left(N-1;\lambda\right)\nonumber \\
 & =-\lambda\left(a_{0}\left(N;\lambda\right)+\left(N+2\lambda\right)a_{0}\left(N-1;\lambda\right)\right)\nonumber \\
 & \relphantom =+\left(N+2\lambda\right)\left(N+2\lambda-1\right)\left\{ -\lambda a_{0}\left(N-2;\lambda\right)+\left(N+2\lambda-2\right)a_{1}\left(N-2;\lambda\right)\right\} \nonumber \\
 & =-\lambda\left\{ a_{0}\left(N;\lambda\right)+\left(N+2\lambda\right)a_{0}\left(N-1;\lambda\right)+\left(N+2\lambda\right)\left(N+2\lambda-1\right)a_{0}\left(N-2;\lambda\right)\right\} \nonumber \\
 & \relphantom =+\left(N+2\lambda\right)\left(N+2\lambda-1\right)\left(N+2\lambda-2\right)a_{1}\left(N-2;\lambda\right)\nonumber \\
 & \vdots\nonumber \\
 & =-\lambda\sum_{i=0}^{N-1}\left(N+2\lambda\right)_{i}a_{0}\left(N-i;\lambda\right)+\left(N+2\lambda\right)_{N}a_{1}\left(1;\lambda\right)\nonumber \\
 & =-\lambda\sum_{i=0}^{N}\left(N+2\lambda\right)_{i}a_{0}\left(N-i;\lambda\right),\nonumber 
\end{align}
\begin{align*}
&\relphantom{=}a_{2}\left(N+1;\lambda\right)\\
& =-2\lambda a_{1}\left(N;\lambda\right)+\left(N+3\lambda\right)a_{2}\left(N;\lambda\right)\\
 & =-2\lambda a_{1}\left(N;\lambda\right)+\left(N+3\lambda\right)\left\{ -2\lambda a_{1}\left(N-1;\lambda\right)+\left(N+3\lambda-1\right)a_{2}\left(N-1;\lambda\right)\right\} \\
 & =-2\lambda\left\{ a_{1}\left(N;\lambda\right)+\left(N+3\lambda\right)a_{1}\left(N-1;\lambda\right)\right\} \\
 & \relphantom =+\left(N+3\lambda\right)\left(N+3\lambda-1\right)\left\{ -2\lambda a_{1}\left(N-2;\lambda\right)+\left(N+3\lambda-2\right)a_{2}\left(N-2;\lambda\right)\right\} \\
 & =-2\lambda\left\{ a_{1}\left(N;\lambda\right)+\left(N+3\lambda\right)a_{1}\left(N-1;\lambda\right)+\left(N+3\lambda\right)\left(N+3\lambda-1\right)a_{1}\left(N-2;\lambda\right)\right\} \\
 & \relphantom =+\left(N+3\lambda\right)\left(N+3\lambda-1\right)\left(N+3\lambda-2\right)a_{2}\left(N-2;\lambda\right)\\
 & \vdots\\
 & =-2\lambda\sum_{i=0}^{N-2}\left(N+3\lambda\right)_{i}a_{1}\left(N-i;\lambda\right)+\left(N+3\lambda\right)_{N-1}a_{2}\left(2;\lambda\right)\\
 & =-2\lambda\sum_{i=0}^{N-1}\left(N+3\lambda\right)_{i}a_{1}\left(N-i;\lambda\right),
\end{align*}
and 
\begin{align}
&\relphantom{=}a_{3}\left(N+1;\lambda\right) \label{eq:21}\\
& =-3\lambda a_{2}\left(N;\lambda\right)+\left(N+4\lambda\right)a_{3}\left(N;\lambda\right)\nonumber\\
 & =-3\lambda \left\{a_{2}\left(N;\lambda\right)+\left(N+4\lambda\right) a_{2}\left(N-1;\lambda\right)\right\}\nonumber \\
 & \relphantom =+\left(N+4\lambda\right)\left(N+4\lambda-1\right)\left\{ -3\lambda a_{2}\left(N-2;\lambda\right)+\left(N+4\lambda-2\right)a_{3}\left(N-2;\lambda\right)\right\} \nonumber \\
 & =-3\lambda\left\{ a_{2}\left(N;\lambda\right)+\left(N+4\lambda\right)a_{2}\left(N-1;\lambda\right)+\left(N+4\lambda\right)\left(N+4\lambda-1\right)a_{2}\left(N-2;\lambda\right)\right\} \nonumber \\
 & \relphantom =+\left(N+4\lambda\right)\left(N+4\lambda-1\right)\left(N+4\lambda-2\right)a_{3}\left(N-2;\lambda\right)\nonumber \\
 & \vdots\\
 & =-3\lambda\sum_{i=0}^{N-3}\left(N+4\lambda\right)_{i}a_{2}\left(N-i;\lambda\right)+\left(N+4\lambda\right)_{N-2}a_{3}\left(3;\lambda\right)\\
 & =-3\lambda\sum_{i=0}^{N-2}\left(N+4\lambda\right)_{i}a_{2}\left(N-i;\lambda\right).
\end{align}

Proceeding in this way, we get 
\begin{equation}
a_{k}\left(N+1;\lambda\right)=-k\lambda\sum_{i_{1}=0}^{N-k+1}\left(N+\left(k+1\right)\lambda\right)_{i_{1}}a_{k-1}\left(N-i_{1};\lambda\right),\label{eq:22}
\end{equation}
where $1\le k\le N$. 

Therefore, we obtain the following theorem.
\begin{thm}
\label{thm:1}  We have the following recurrence relations:
\begin{enumerate}
\item[(i)]  $a_{0}\left(0;\lambda\right)=\frac{1}{\lambda}$, $a_{0}\left(1;\lambda\right)=1$,
$a_{1}\left(1;\lambda\right)=-1$,
\item[(ii)]  $a_{0}\left(N+1;\lambda\right)=\left(N+\lambda\right)_{N}$, $a_{N+1}\left(N+1;\lambda\right)=\left(-1\right)^{N+1}\lambda^{N}\left(N+1\right)!$,
\item[(iii)]  $a_{k}\left(N+1;\lambda\right)=-k\lambda\sum_{i_{1}=0}^{N-k+1}\left(N+\left(k+1\right)\lambda\right)_{i_{1}}a_{k-1}\left(N-i_{1};\lambda\right)$, 
\end{enumerate}
\end{thm}
for $1\le k\le N$.

Now, we observe that 
\begin{align}
a_{1}\left(N+1;\lambda\right) & =-\lambda\sum_{i_{1}=0}^{N}\left(N+2\lambda\right)_{i_{1}}a_{0}\left(N-i_{1};\lambda\right)\label{eq:23}\\
 & =-\lambda\sum_{i_{1}=0}^{N}\left(N+2\lambda\right)_{i_{1}}\left(N+\lambda-i_{1}-1\right)_{N-i_{1}-1},\nonumber 
\end{align}
\begin{align}
&\relphantom{=}a_{2}\left(N+1;\lambda\right) \label{eq:24}\\
& =-2\lambda\sum_{i_{2}=0}^{N-1}\left(N+3\lambda\right)_{i_{2}}a_{1}\left(N-i_{2};\lambda\right)\\
 & =\left(-1\right)^{2}2!\lambda^{2}\sum_{i_{2}=0}^{N-1}\sum_{i_{1}=0}^{N-i_{2}-1}\left(N+3\lambda\right)_{i_{2}}\left(N+2\lambda-i_{2}-1\right)_{i_{1}}\nonumber\\
&\relphantom{=}\times\left(N+\lambda-i_{2}-i_{1}-2\right)_{N-i_{2}-i_{1}-2},\nonumber 
\end{align}
and 
\begin{align}
&a_{3}\left(N+1;\lambda\right) \label{eq:25}\\
& =-3\lambda\sum_{i_{3}=0}^{N-2}\left(N+4\lambda\right)_{i_{3}}a_{2}\left(N-i_{3};\lambda\right)\nonumber\\
 & =\left(-1\right)^{3}3!\lambda^{3}\sum_{i_{3}=0}^{N-2}\sum_{i_{2}=0}^{N-i_{3}-2}\sum_{i_{1}=0}^{N-i_{3}-i_{2}-2}\left(N+4\lambda\right)_{i_{3}}\left(N+3\lambda-i_{3}-1\right)_{i_{2}}\nonumber\\
&\relphantom{=}\times\left(N+2\lambda-i_{3}-i_{2}-2\right)_{i_{1}}\nonumber \\
 & \relphantom =\times\left(N+\lambda-i_{3}-i_{2}-i_{1}-3\right)_{N-i_{3}-i_{2}-i_{1}-3}.\nonumber 
\end{align}

Continuing this process, we have 
\begin{align}
 & \relphantom =a_{j}\left(N+1;\lambda\right)\label{eq:26}\\
 & =\left(-1\right)^{j}j!\lambda^{j}\nonumber\\
&\relphantom{=}\times \sum_{i_{j}=0}^{N-j+1}\sum_{i_{j-1}=0}^{N-j+1-i_{j}}\cdots\sum_{i_{1}=0}^{N-j+1-i_{j}-\cdots-i_{2}}\left(N+\left(j+1\right)\lambda\right)_{i_{j}}\left(N+j\lambda-i_{j}-1\right)_{j-1}\nonumber \\
 & \relphantom =\times\cdots\times\left(N+2\lambda-i_{j}-\cdots-i_{2}-\left(j-1\right)\right)_{i_{1}}\nonumber\\
&\relphantom{=}\times\left(N+\lambda-i_{j}-\cdots-i_{1}-j\right)_{N-i_{j}-\cdots-i_{1}-j},\nonumber 
\end{align}
where $1\le j\le N$. 

From (\ref{eq:26}), we note that the matrix $\left(a_{i}\left(j;\lambda\right)\right)_{0\le i,j\le N}$
is given by 
\begin{equation}
 \begin{tikzpicture}[baseline=(current  bounding  box.west)]
  \matrix (mymatrix) [matrix of math nodes,left delimiter={[},right
delimiter={]}]
  {
    \frac{1}{\lambda}  & 1   & (1+\lambda) & (2+\lambda)_2 & \cdots & (N+\lambda-1)_{N-1} \\
      & -1   &  &  &    \\
     &  &  (-1)^2 \lambda 2! &  &     \\
     &  & & (-1)^3 \lambda^2 3! & &  \\
    &  &  & &\ddots &  \\
    &  &  & & & (-1)^N \lambda^{N-1}N! \\
  };
\node[xshift=-105pt,yshift=68pt] {$0$};
\node[xshift=-88pt,yshift=68pt] {$1$};
\node[xshift=-57pt,yshift=68pt] {$2$};
\node[xshift=-12pt,yshift=68pt] {$3$};
\node[xshift=70pt,yshift=68pt] {$N$};
\node[xshift=-135pt,yshift=43pt] {$0$};
\node[xshift=-135pt,yshift=28pt] {$1$};
\node[xshift=-135pt,yshift=11pt] {$2$};
\node[xshift=-135pt,yshift=-10pt] {$3$};
\node[xshift=-135pt,yshift=-45pt] {$N$};
\node[xshift=-80pt,yshift=-30pt] {\LARGE$0$};
\end{tikzpicture}\label{eq:27}
\end{equation}

Therefore, by Theorem 1, (\ref{eq:7}), and (\ref{eq:26}), we obtain the following theorem.
\begin{thm}
\label{thm:2} The nonlinear differential equations 
\[
F^{\left(N\right)}=\frac{\left(-1\right)^{N}\lambda}{\left(1+t\right)^{N}}\sum_{i=1}^{N+1}a_{i-1}\left(N;\lambda\right)F^{i},\quad\left(N\in\mathbb{N}\right),
\]
have a solution $F=F\left(t,\lambda\right)=\frac{1}{\left(1+t\right)^{\lambda}+1}$,\\
where $a_{0}\left(N;\lambda\right)=\left(N+\lambda-1\right)_{N-1}$, $a_{N}\left(N;\lambda\right)=\left(-1\right)^{N}\lambda^{N-1}N!$,
\begin{align*}
 & a_{j}\left(N;\lambda\right)\\
 & =\left(-1\right)^{j}j!\lambda^{j}\sum_{i_{j}=0}^{N-j}\sum_{i_{j-1}=0}^{N-j-i_{j}}\cdots\sum_{i_{1}=0}^{N-j-i_{j}-\cdots-i_{2}}\left(N+\left(j+1\right)\lambda-1\right)_{i_{j}}\\
 & \relphantom =\times\left(N+j\lambda-\lambda_{j}-2\right)_{i_{j-1}}\cdots\left(N+2\lambda-i_{j}-\cdots-i_{2}-j\right)_{i_{1}}\\
 & \relphantom =\times\left(N+\lambda-i_{j}-\cdots-i_{1}-j-1\right)_{N-i_{j}-\cdots-i_{1}-j-1}, \quad \left(1\le j\le N-1\right).
\end{align*}

\end{thm}
Recall that the Boole numbers, $Bl_{k}\left(\lambda\right)$, $\left(k\ge0\right)$,
are given by the generating function  
\begin{equation}
\frac{1}{\left(1+t\right)^{\lambda}+1}=\sum_{k=0}^{\infty}Bl_{k}\left(\lambda\right)\frac{t^{k}}{k!}.\label{eq:28}
\end{equation}

Thus, by (\ref{eq:28}), we get 
\begin{align*}
 & F^{\left(N\right)}\\
 & =\left(\frac{d}{dt}\right)^{N}F\left(t,\lambda\right)\\
 & =\left(\frac{d}{dt}\right)^{N}\left(\frac{1}{\left(1+t\right)^{\lambda}+1}\right)\\
 & =\sum_{k=N}^{\infty}Bl_{k}\left(\lambda\right)\left(k\right)_{N}\frac{t^{k-N}}{k!}\\
 & =\sum_{k=0}^{\infty}Bl_{k+N}\left(\lambda\right)\frac{\left(k+N\right)_{N}}{\left(k+N\right)!}t^{k}\\
 & =\sum_{k=0}^{\infty}Bl_{k+N}\left(\lambda\right)\frac{t^{k}}{k!},\quad\left(N\in\mathbb{N}\right).
\end{align*}

From (\ref{eq:2}), Theorem \ref{thm:2} and (\ref{eq:28}), we have
\begin{align}
&\relphantom{=}\sum_{k=0}^{\infty}Bl_{k+N}\left(\lambda\right)\frac{t^{k}}{k!} \label{eq:29}\\
& =F^{\left(N\right)}\nonumber\\
 & =\frac{\left(-1\right)^{N}\lambda}{\left(1+t\right)^{N}}\sum_{i=1}^{N+1}a_{i-1}\left(N;\lambda\right)\left(\frac{1}{\left(1+t\right)^{\lambda}+1}\right)^{i}\nonumber \\
 & =\left(-1\right)^{N}\lambda\left(1+t\right)^{-N}\sum_{i=1}^{N+1}a_{i-1}\left(N;\lambda\right)\left(\frac{1}{\left(1+t\right)^{\lambda}+1}\right)^{i}\nonumber \\
 & =\left(-1\right)^{N}\lambda\left(\sum_{l=0}^{\infty}\left(-1\right)^{l}\left(N+l-1\right)_{l}\frac{t^{l}}{l!}\right)\left(\sum_{i=1}^{N+1}a_{i-1}\left(N;\lambda\right)\sum_{m=0}^{\infty}Bl_{m}^{\left(i\right)}\left(\lambda\right)\frac{t^{m}}{m!}\right)\nonumber \\
 & =\left(-1\right)^{N}\lambda\sum_{i=1}^{N+1}a_{i-1}\left(N;\lambda\right)\left(\sum_{l=0}^{\infty}\left(-1\right)^{l}\left(N+l-1\right)_{l}\frac{t^{l}}{l!}\right)\left(\sum_{m=0}^{\infty}Bl_{m}^{\left(i\right)}\left(\lambda\right)\frac{t^{m}}{m!}\right)\nonumber \\
 & =\left(-1\right)^{N}\lambda\sum_{i=1}^{N+1}a_{i-1}\left(N;\lambda\right)\left(\sum_{k=0}^{\infty}\sum_{l=0}^{k}\binom{k}{l}\left(-1\right)^{l}\left(N+l-1\right)_{l}Bl_{k-l}^{\left(i\right)}\left(\lambda\right)\right)\frac{t^{k}}{k!}\nonumber \\
 & =\sum_{k=0}^{\infty}\left\{ \left(-1\right)^{N}\lambda\sum_{i=1}^{N+1}a_{i-1}\left(N;\lambda\right)\sum_{l=0}^{k}\binom{k}{l}\left(-1\right)^{l}\left(N+l-1\right)_{l}Bl_{k-l}^{\left(i\right)}\left(\lambda\right)\right\} \frac{t^{k}}{k!},\nonumber 
\end{align}
where $N\in\mathbb{N}$.

By comparing the coefficients on both sides of (\ref{eq:29}),
we obtain the following theorem.
\begin{thm}
\label{thm:3} For $N\in\mathbb{N}$ and $k\in\mathbb{N}\cup\left\{ 0\right\} $,
we have 
\[
Bl_{k+N}\left(\lambda\right)=\left(-1\right)^{N}\lambda\sum_{i=1}^{N+1}a_{i-1}\left(N;\lambda\right)\sum_{k=0}^{k}\binom{k}{l}\left(-1\right)^{l}\left(N+l-1\right)_{l}Bl_{k-l}^{\left(i\right)}\left(\lambda\right).
\]

\end{thm}
By replacing $t$ by $e^{t}-1$ in (\ref{eq:1}), we get 
\begin{align}
\frac{1}{2}\left(\frac{2}{e^{\lambda t}+1}\right) & =\sum_{k=0}^{\infty}Bl_{k}\left(\lambda\right)\frac{1}{k!}\left(e^{t}-1\right)^{k}\label{eq:30}\\
 & =\sum_{n=0}^{\infty}\left(\sum_{k=0}^{n}Bl_{k}\left(\lambda\right)S_{2}\left(n,k\right)\right)\frac{t^{n}}{n!},\nonumber 
\end{align}
where $S_{2}\left(n,k\right)$ are the Stirling numbers of the second
kind.

As is well known, Euler numbers are given by the generating function

\begin{equation}
\left(\frac{2}{e^{t}+1}\right)=\sum_{n=0}^{\infty}E_{n}\frac{t^{n}}{n!},\quad\left(\text{see \cite{key-7}}\right).\label{eq:31}
\end{equation}

From (\ref{eq:30}) and (\ref{eq:31}), we have 
\begin{equation}
2^{-1}\lambda^{n}E_{n}=\sum_{k=0}^{n}Bl_{k}\left(\lambda\right)S_{2}\left(n,k\right),\quad\left(n\ge0\right).\label{eq:32}
\end{equation}

It is well known that the higher-order Euler numbers are also defined
by the generating function 
\begin{equation}
\left(\frac{2}{e^{t}+1}\right)^{r}=\sum_{n=0}^{\infty}E_{n}^{\left(r\right)}\frac{t^{n}}{n!},\quad\left(\text{see \cite{key-19}}\right).\label{eq:33}
\end{equation}

Now, we observe that 
\begin{align}
\left(\frac{1}{e^{\lambda t}+1}\right)^{i} & =\left(\frac{1}{\left(e^{t}-1+1\right)^{\lambda}+1}\right)^{i}\label{eq:34}\\
 & =\sum_{k=0}^{\infty}Bl_{k}^{\left(i\right)}\left(\lambda\right)\frac{1}{k!}\left(e^{t}-1\right)^{k}\nonumber \\
 & =\sum_{n=0}^{\infty}\left(\sum_{k=0}^{n}Bl_{k}^{\left(i\right)}\left(\lambda\right)S_{2}\left(n,k\right)\right)\frac{t^{n}}{n!}.\nonumber 
\end{align}

Thus, by (\ref{eq:33}) and (\ref{eq:34}), we get 
\[
2^{-i}\lambda^{n}E_{n}^{\left(i\right)}=\sum_{k=0}^{n}Bl_{k}^{\left(i\right)}\left(\lambda\right)S_{2}\left(n,k\right),\quad\left(n\ge0,i\in\mathbb{N}\right).
\]

From (\ref{eq:1}) and (\ref{eq:31}), we note that 
\begin{align}
2\sum_{n=0}^{\infty}Bl_{n}\left(\lambda\right)\frac{t^{n}}{n!} & =\frac{2}{\left(1+t\right)^{\lambda}+1}\label{eq:35}\\
 & =\frac{2}{e^{\lambda\log\left(1+t\right)}+1}\nonumber \\
 & =\sum_{k=0}^{\infty}E_{k}\frac{\lambda^{k}}{k!}\left(\log\left(1+t\right)\right)^{k}\nonumber \\
 & =\sum_{n=0}^{\infty}\left(\sum_{k=0}^{n}E_{k}\lambda^{k}S_{1}\left(n,k\right)\right)\frac{t^{n}}{n!},\nonumber 
\end{align}
where $S_{1}\left(n,k\right)$ are the Stirling numbers of the first
kind. 

Thus, by (\ref{eq:35}), we get 
\begin{equation}
Bl_{n}\left(\lambda\right)=\frac{1}{2}\sum_{k=0}^{n}E_{k}\lambda^{k}S_{1}\left(n,k\right),\quad\left(n\ge0\right).\label{eq:36}
\end{equation}

By (\ref{eq:33}), we easily get 
\begin{align}
\left(\frac{2}{\left(1+t\right)^{\lambda}+1}\right)^{i} & =\left(\frac{2}{e^{\lambda\log\left(1+t\right)}+1}\right)^{i}\label{eq:37}\\
 & =\sum_{k=0}^{\infty}E_{k}^{\left(i\right)}\frac{1}{k!}\lambda^{k}\left(\log\left(1+t\right)\right)^{k}\nonumber \\
 & =\sum_{n=0}^{\infty}\left(\sum_{k=0}^{n}E_{k}^{\left(i\right)}\lambda^{k}S_{1}\left(n,k\right)\right)\frac{t^{n}}{n!},\quad\left(i\in\mathbb{N}\right).\nonumber 
\end{align}

From (\ref{eq:2}) and (\ref{eq:37}), we have 
\begin{equation}
2^{i}Bl_{n}^{\left(i\right)}\left(\lambda\right)=\sum_{k=0}^{n}E_{k}^{\left(i\right)}\lambda^{k}S_{1}\left(n,k\right),\quad\left(n\ge0,i\in\mathbb{N}\right).\label{eq:38}
\end{equation}

Therefore, by Theorem \ref{thm:3}, (\ref{eq:37}), and (\ref{eq:38}),
we obtain the following theorem.
\begin{thm}
\label{thm:4} For $k\in\mathbb{N}\cup\left\{ 0\right\} $ and $N\in\mathbb{N}$,
we have 
\begin{align*}
 & \frac{1}{2}\sum_{n=0}^{k+N}E_{n}\lambda^{n}S_{1}\left(k+N,n\right)\\
 & =\left(-1\right)^{N}\lambda\sum_{i=1}^{N+1}a_{i-1}\left(N;\lambda\right)\sum_{l=0}^{k}\binom{k}{l}\left(-1\right)^{l}\left(N+l-1\right)_{l}\sum_{n=0}^{k-l}2^{-i}E_{n}^{\left(i\right)}\lambda^{n}S_{1}\left(k-l,n\right).
\end{align*}
\end{thm}

\bibliographystyle{amsplain}
\providecommand{\bysame}{\leavevmode\hbox to3em{\hrulefill}\thinspace}
\providecommand{\MR}{\relax\ifhmode\unskip\space\fi MR }
\providecommand{\MRhref}[2]{%
  \href{http://www.ams.org/mathscinet-getitem?mr=#1}{#2}
}
\providecommand{\href}[2]{#2}

\end{document}